\documentclass{article}
\usepackage{graphicx}

\newtheorem{theorem}{Theorem}

\newtheorem{definition}[theorem]{Definition}

\title {Auxetic regions in large deformations of periodic frameworks}
\author{Ciprian S. Borcea and Ileana Streinu}

\begin{document}
\maketitle

\begin{abstract}
In materials science, auxetic behavior refers to lateral widening upon stretching. We investigate the problem of finding domains of auxeticity in global deformation spaces of periodic frameworks. Case studies include planar periodic mechanisms constructed from quadrilaterals with diagonals as periods and other frameworks with two vertex orbits. We relate several geometric and kinematic descriptions.
\end{abstract}


\section{Introduction}
\label{sec:1}

Periodic frameworks are mathematical abstractions which model atom-and-bond structures in crystalline materials or man-made trusses and repetitive articulated systems \cite{BS1,Su}.
Recent advances in additive manufacturing have brought about new possibilities for 
producing complex structures at microscales, giving new impetus to rational design of metamaterials, particularly those based on a periodic organization \cite{BVCvH,CKKW,RJvH}. 

\medskip
Auxetic behavior is a rather counter-intuitive mechanical response of certain
materials and structures which involves lateral widening
under stretching and, in reverse, lateral shrinking under compression. 
This type of response is desirable for various applications, ranging from medical implements to shock-absorbing curtains and has been investigated with increased interest over the last three decades \cite{GGLR,KZ,L,SDC}.
Even so, most studies have addressed only
small deformations of a confined repertory of structural designs and the problem of
large deformations has been barely considered \cite{ZLWRAD}. Auxetic growth implies, in particular, volume increase \cite{BS6} and is therefore limited by certain structural bounds. Thus, for large deformations, auxetic behavior 
can be expected only over restricted regions. The main purpose of the 
present contribution is to offer an introduction to this type of kinematic inquiry,
based on our geometric theory of auxetic deformations and recent fundamental results on the structure and design of auxetic periodic frameworks \cite{BS6,BS9,BS10}.

\medskip
We review in section~\ref{sec:review} some key notions  
on periodic frameworks and auxetics.
We proceed in section~\ref{sec:4bar} with a fairly intuitive scenario in dimension two, by associating a 
periodic framework to a given quadrilateral, which is repeatedly translated by its diagonal vectors, as shown below in Figure~\ref{fig:4barQ}. 
As a planar linkage, a quadrilateral is a four-bar mechanism and our opening discussion is just an interpretation of the classical kinematics of this mechanism in the periodic setting. Yet, even this elementary situation proves the value of a strictly geometric approach to auxetics and suggests useful algebraic
parametrizations based on spectrahedra. This point of view is pursued in   subsequent sections and leads to a unified perspective for identifying auxetic regions in global deformations of periodic frameworks with
two vertex orbits.

\section{Prerequisites}
\label{sec:review}

We begin with a concise review of notions which are necessary for a mathematical formulation of our problem.

\medskip \noindent
{\bf Periodic frameworks and deformation spaces.} The principal reference is \cite{BS1}.

\begin{definition}
	A {\em $d$-periodic graph} is a pair $(G,\Gamma)$, where $G=(V,E)$ is a simple infinite graph with vertices $V$, edges
	$E$ and finite degree at every vertex, and $\Gamma \subset Aut(G)$ is a free Abelian group of automorphisms
	which has rank $d$, acts without fixed points and has a finite number of vertex (and hence, also edge) orbits.
\end{definition}
 
\medskip \noindent
We assume $G$ to be connected. The group $\Gamma$ is isomorphic to $ {Z}^d$ and is called the {\em periodicity group}  of the periodic graph $G$. Vertices which are
equivalent under $\Gamma$ form a vertex orbit and similarly for edges. 

\begin{definition}
	A (periodic) {\em placement} of a $d$-periodic graph $(G,\Gamma)$ in 
	$  {R}^d$ is defined by two functions:
	$p:V\rightarrow   {R}^d$ and  
	$\pi: \Gamma \hookrightarrow {\cal T}(  {R}^d)$, 
	where $p$ assigns points in $  {R}^d$ to the vertices $V$ of $G$ and $\pi$ is a faithful representation of the periodicity group $\Gamma$, that is, an injective homomorphism of $\Gamma$ into the group ${\cal T}( {R}^d)$ of translations in the Euclidean space $  {R}^d$, with $\pi(\Gamma)$ being a lattice of rank $d$. These two functions must satisfy the natural compatibility condition $p(\gamma v)=\pi(\gamma)(p(v))$. The translation group ${\cal T}( {R}^d)$ can be identified with the additive group of vectors in $ {R}^d$.
\end{definition}

\medskip \noindent
In a framework, edges are seen as segments between the corresponding vertices. Mechanical interpretations consider edges as rigid bars and vertices as (spherical) joints, hence the name bar-and-joint framework.

\begin{definition}
	Given a $d$-periodic framework ${\cal F}=(G,\Gamma, p,\pi)$, the collection of all periodic placements of $(G,\Gamma)$ in ${R}^d$ which maintain the lengths of all edges is called the {\em realization space} of the framework. After factoring out equivalence under Euclidean isometries, we obtain the {\em configuration space} of the framework (with the quotient topology). The {\em deformation space} is the connected component of the configuration space which contains the initial framework.
\end{definition}

\medskip \noindent
It is important to notice that the representation of the periodicity group $\Gamma$
by a lattice of translations may well vary when the framework deforms.

\medskip \noindent
{\bf Auxetic deformations.} The motivations for the following definition of auxetic paths are detailed in \cite{BS6}. Let $(G,\Gamma, p_{\tau},\pi_{\tau}), \tau\in (-\epsilon,\epsilon)$ be a one-parameter deformation of the periodic
framework $(G,\Gamma, p, \pi)$, where $p=p_0$ and $\pi=\pi_0$. Suppose we have
chosen an independent set of generators $\gamma_i$, $i=1,...,d$ for the periodicity group $\Gamma$. Then, we have at every moment $\tau$ a lattice basis $\pi_{\tau}(\gamma_i)$, which we write as a $d\times d$ matrix (with column vectors) $\Lambda_{\tau}$. The Gram matrix of the lattice basis is therefore
$\omega_{\tau}=\Lambda^t_{\tau}\Lambda_{\tau}$.

\begin{definition}\label{A-path}
	A one-parameter deformation  $(G,\Gamma, p_{\tau},\pi_{\tau}), \tau\in (-\epsilon,\epsilon)$ is called an {\em auxetic path}, or simply {\em auxetic}, when the curve of Gram matrices $\omega_{\tau}$ has all its velocity vectors
	  in the cone of positive semidefinite symmetric $d\times d$ matrices. When all velocity vectors are in the positive definite cone, the deformation is called {\em strictly auxetic.}
\end{definition}

\medskip \noindent
In short, for periodic frameworks, auxetic behavior is expressed
through the evolution of the Gram matrix of periodicity generators. 
A number of recent results derived from this geometric approach to auxetics can be found in \cite{BS9,BS10}. The present inquiry will be oriented towards questions involving large deformations, that go beyond a small neighborhood of a given
initial configuration.

\section{Four-bar mechanisms and associated periodic frameworks}
\label{sec:4bar}
 
\begin{figure}[ht]
\centering
 {\includegraphics[width=0.51\textwidth]{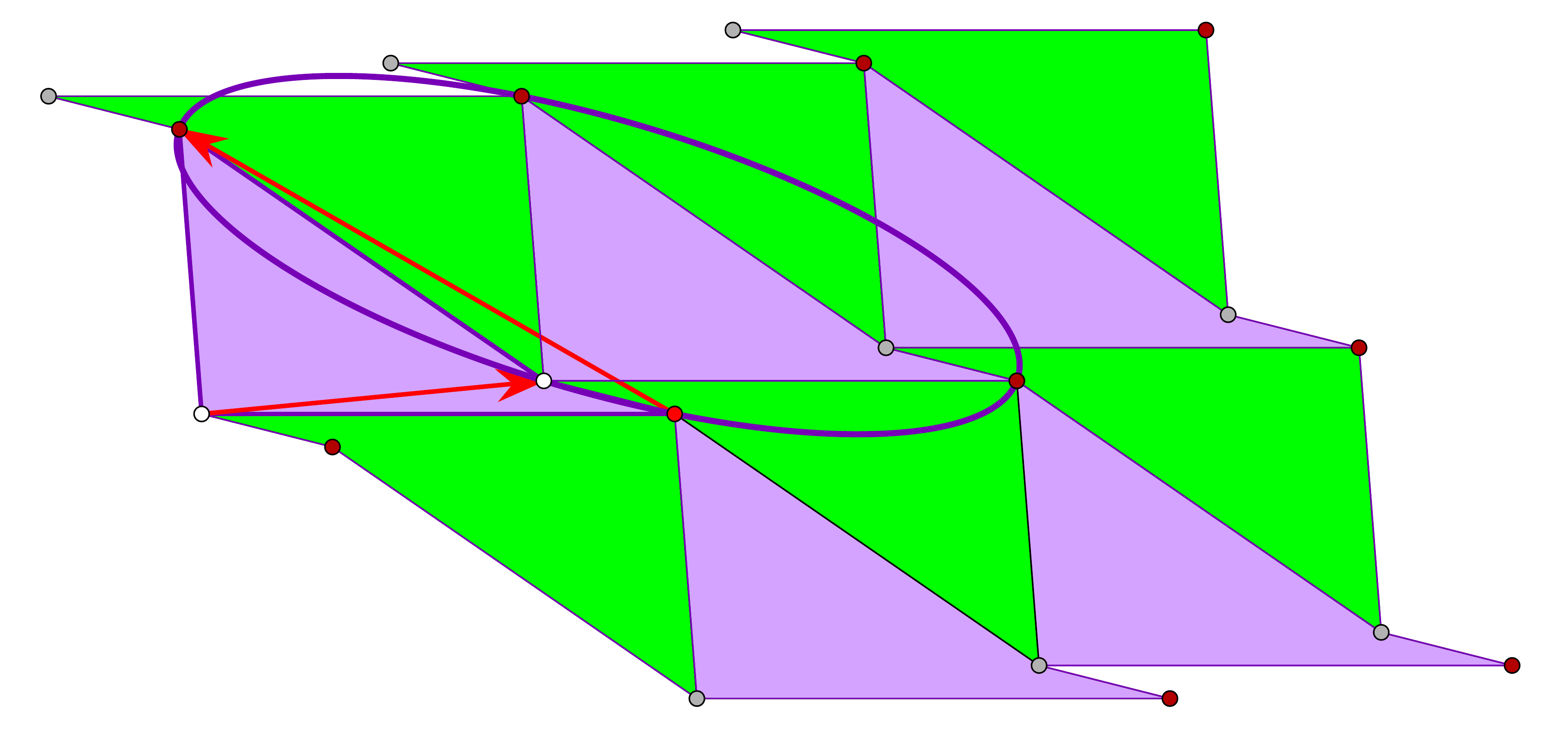}}
 {\includegraphics[width=0.48\textwidth]{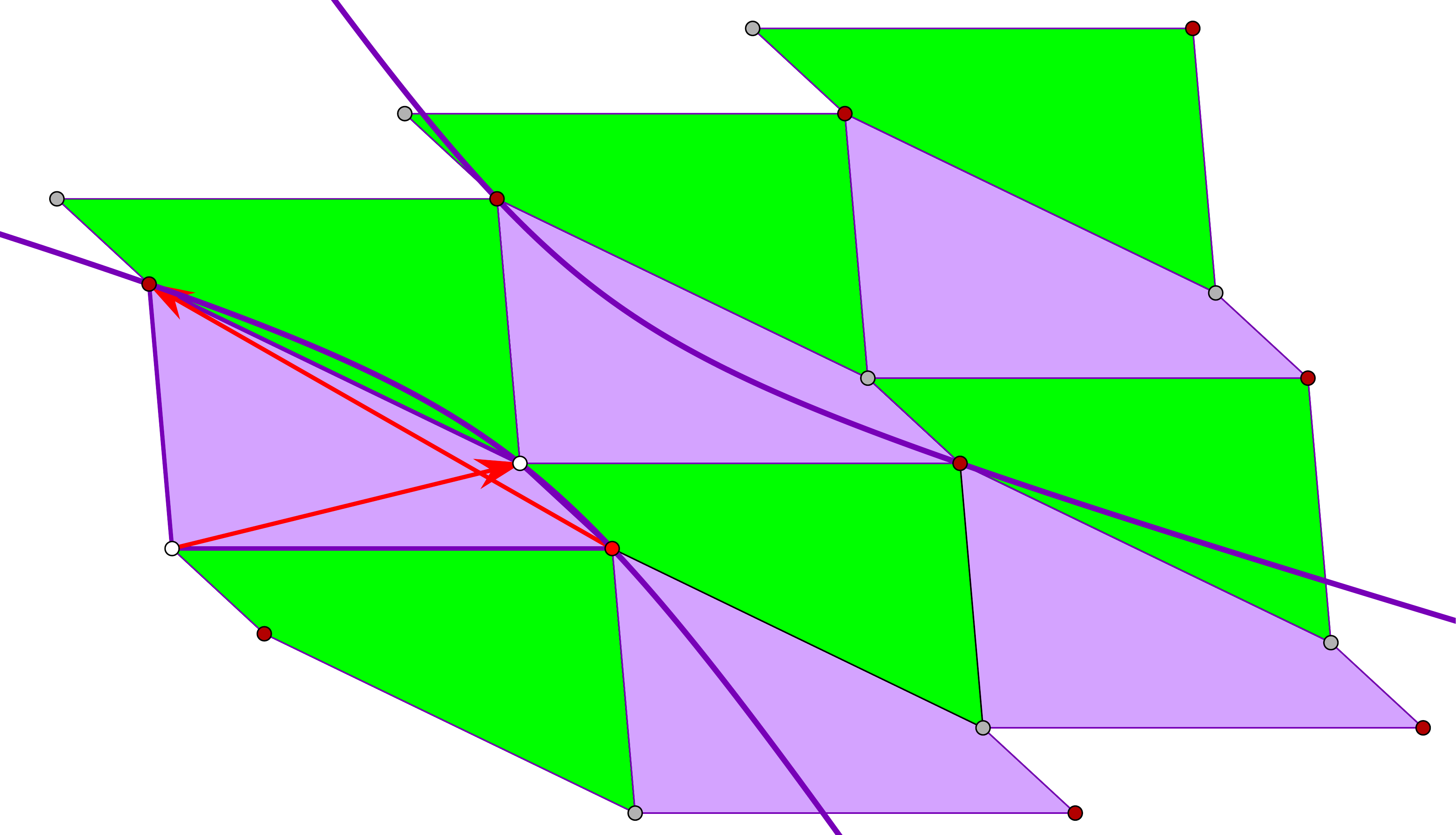}}
 \caption{(Left) A periodic framework can be associated to a quadrilateral by adopting the two diagonals as generators of the periodicity lattice. A four bar mechanism is converted into a planar periodic framework with one degree of freedom. When the quadrilateral is in a pseudotriangle configuration, the result is a periodic pseudotriangulation. The local deformation is expansive, hence auxetic, and this property is confirmed by the ellipse running through a vertex and the endpoints of the four bars emanating from it. (Right) Auxetic behavior ends when the generating quadrilateral becomes convex or self-intersecting. This is reflected in the fact that the conic through a vertex and the endpoints of the four bars emanating from it becomes a hyperbola.}
 \label{fig:4barQ}        
 \end{figure}

We open our investigation with a simple type of planar periodic framework which offers an intuitive setting for the abstract concepts reviewed above. We start with
a quadrilateral $ABCD$ and construct an associated periodic framework by adopting the two diagonal vectors $\overrightarrow{AC}$ and $\overrightarrow{BD}$ as 
periodicity generators. Repeated translations of the quadrilateral by these vectors will articulate into a periodic framework when shared vertices are identified as single joints. The process is illustrated in Fig.~\ref{fig:4barQ}, 
where the two diagonal vectors are shown as red arrows. 
The formal construction is described in a more general setting in \cite{BS9}. 
When the quadrilateral is imagined as a four-bar mechanism (and as long as the two
diagonal vectors remain independent), this association converts a one-degree-of-freedom linkage into a one-degree-of-freedom periodic framework.
 
The deformation space of a four-bar mechanism is well studied in classical kinematics and our task in the periodic setting is reduced to recognizing which intervals correspond to auxetic (or reversed-auxetic)
behavior.

\medskip
We assume a generic quadrilateral. By Grashof's rule, the oriented configuration space is topologically either
 (i) a single loop, if the longest edge plus the shortest edge are more than the other two edges together, or
 (ii) two loops, if the longest edge plus the shortest edge are less than the rest.
 Parallel diagonals occur only in the former case.
Equivalence under reflections identifies the two loops in one case  and gives an 
involution without fixed points of the single loop in the other case (with matching of the two configurations with parallel diagonals).
Thus, in both cases, the deformation space of the linkage (considered up to planar isometries) is topologically a single loop (i.e. a circle). For the deformation space of the associated periodic framework, case (i) requires the exclusion of one point for the case when the diagonals become parallel.

\medskip
For brevity, we expand here on case (ii), which is illustrated in Fig.~\ref{fig:4barQ}.
As mentioned in the introduction,
an auxetic deformation will increase the area of a unit cell (fundamental domain under periodicity) and auxetic intervals in our deformation loop cannot contain
relative minima or maxima of the area function. A precise determination of these
intervals can be obtained by direct computation, but we present a more advanced
point of view, based on pseudotriangulations \cite{BS4,BS5,S}.

\medskip
A pseudotriangle is a simple planar polygon with exactly three interior angles less than $\pi$. In a quadrilateral linkage, the lengths of the two diagonals vary in the same way (i.e. both increase or both decrease) if and only if the quadrilateral is in a pseudotriangle configuration \cite{S}.
Periodic pseudotriangulations are introduced and treated at large in \cite{BS5} and their role in understanding expansive planar periodic deformations is described in \cite{BS4}. Expansive deformations require all distances between pairs of vertices to increase or stay the same. Periodic pseudotriangulations deform expansively and expansive implies auxetic \cite{BS5,BS6}. Assembling these facts, we obtain the following geometric characterization.

\begin{theorem}\label{pseudo}
The periodic framework associated to a quadrilateral has a local auxetic deformation if and only if the quadrilateral is in a pseudotriangle configuration.
\end{theorem}

\noindent
If we use a general structural result established in \cite{BS10}, we can state an equivalent characterization.

\begin{theorem}\label{ellipse}
The periodic framework associated to a quadrilateral has a local auxetic deformation if and only if the conic passing through a chosen vertex and the endpoints of the four bars emanating from it is an ellipse.
\end{theorem}

\begin{figure}[ht]
 \vspace{-2pt}
\centering
{\includegraphics[width=0.23\textwidth]{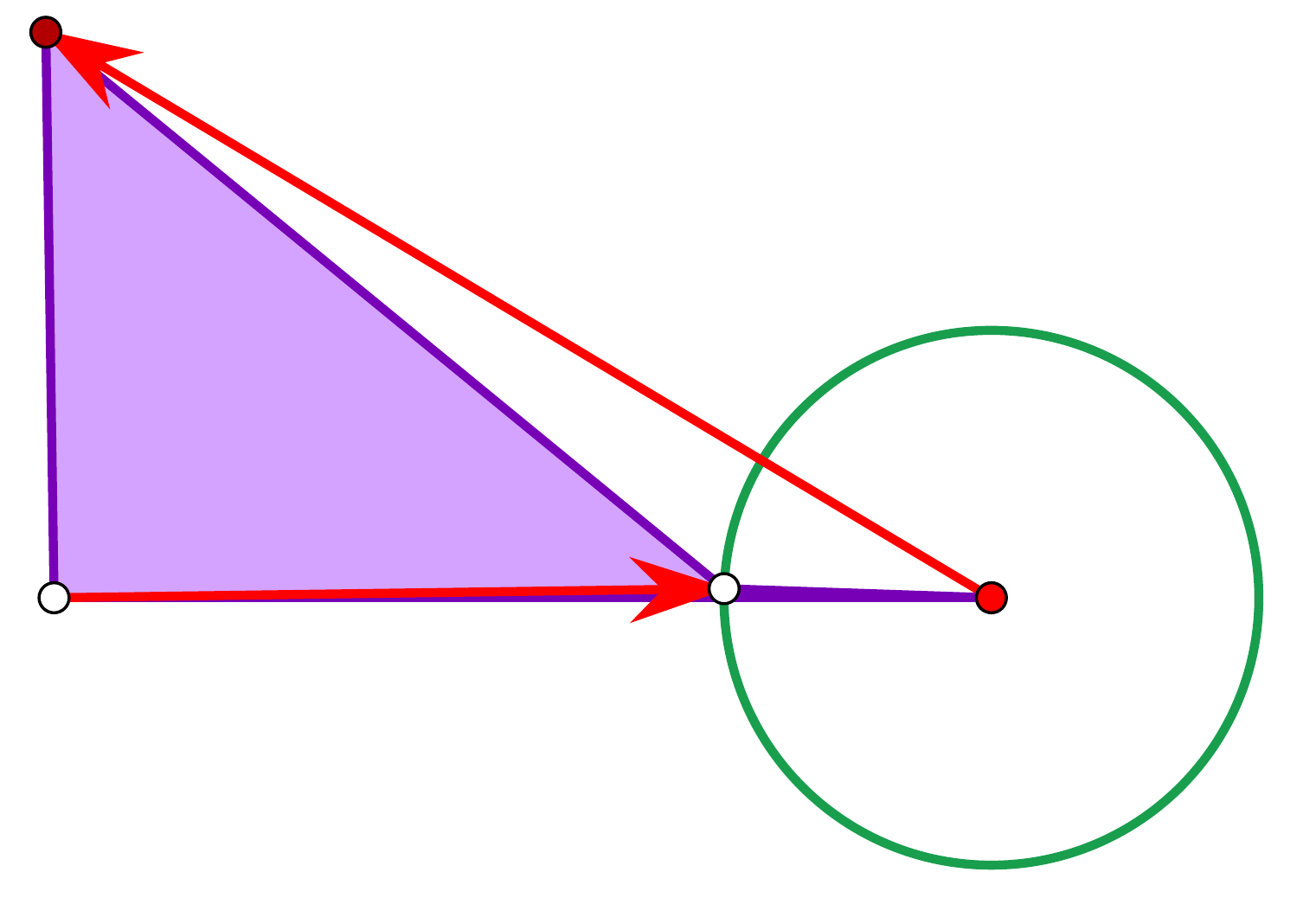}}
\hspace{3pt}
{\includegraphics[width=0.23\textwidth]{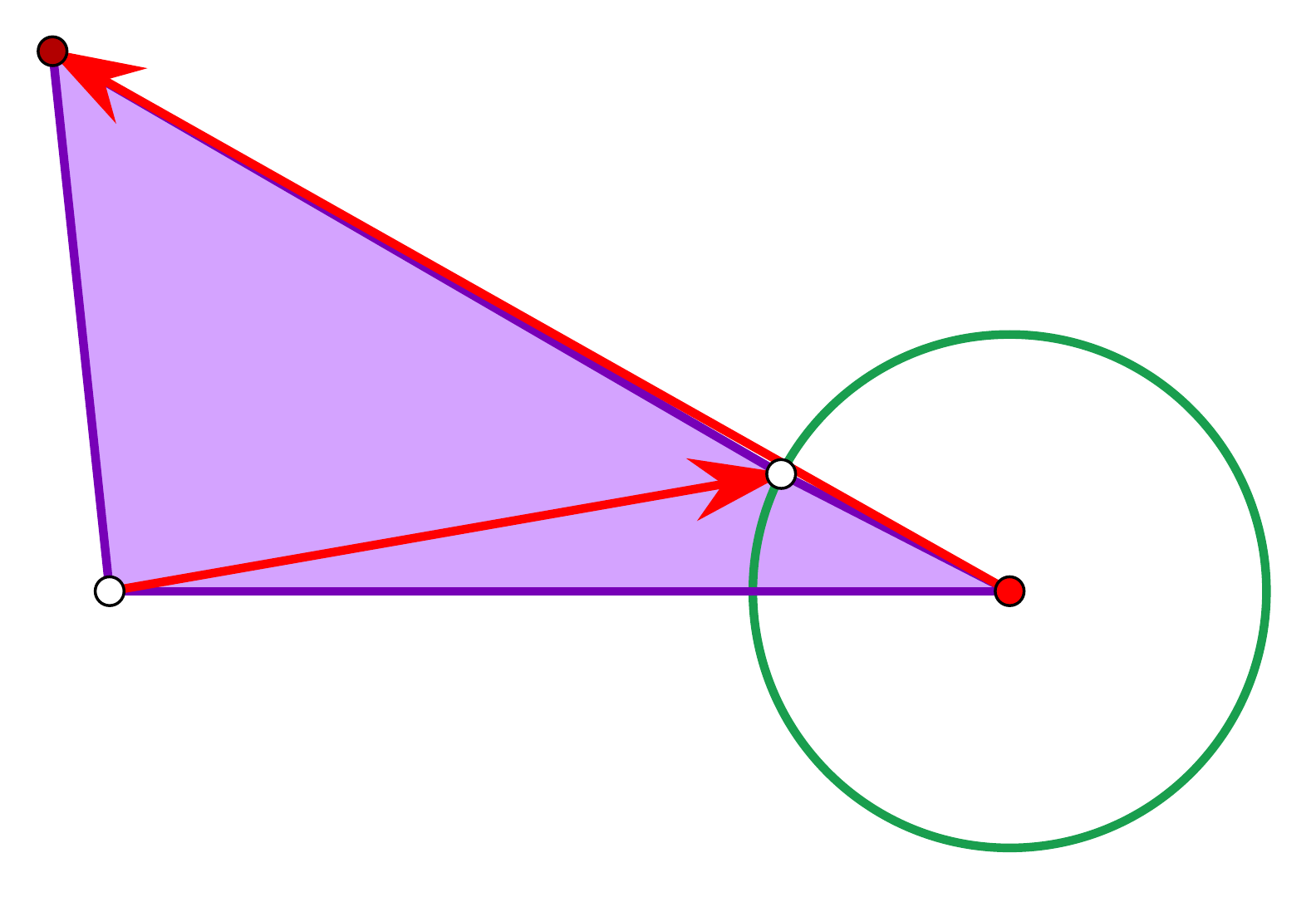}}
\hspace{3pt}
{\includegraphics[width=0.22\textwidth]{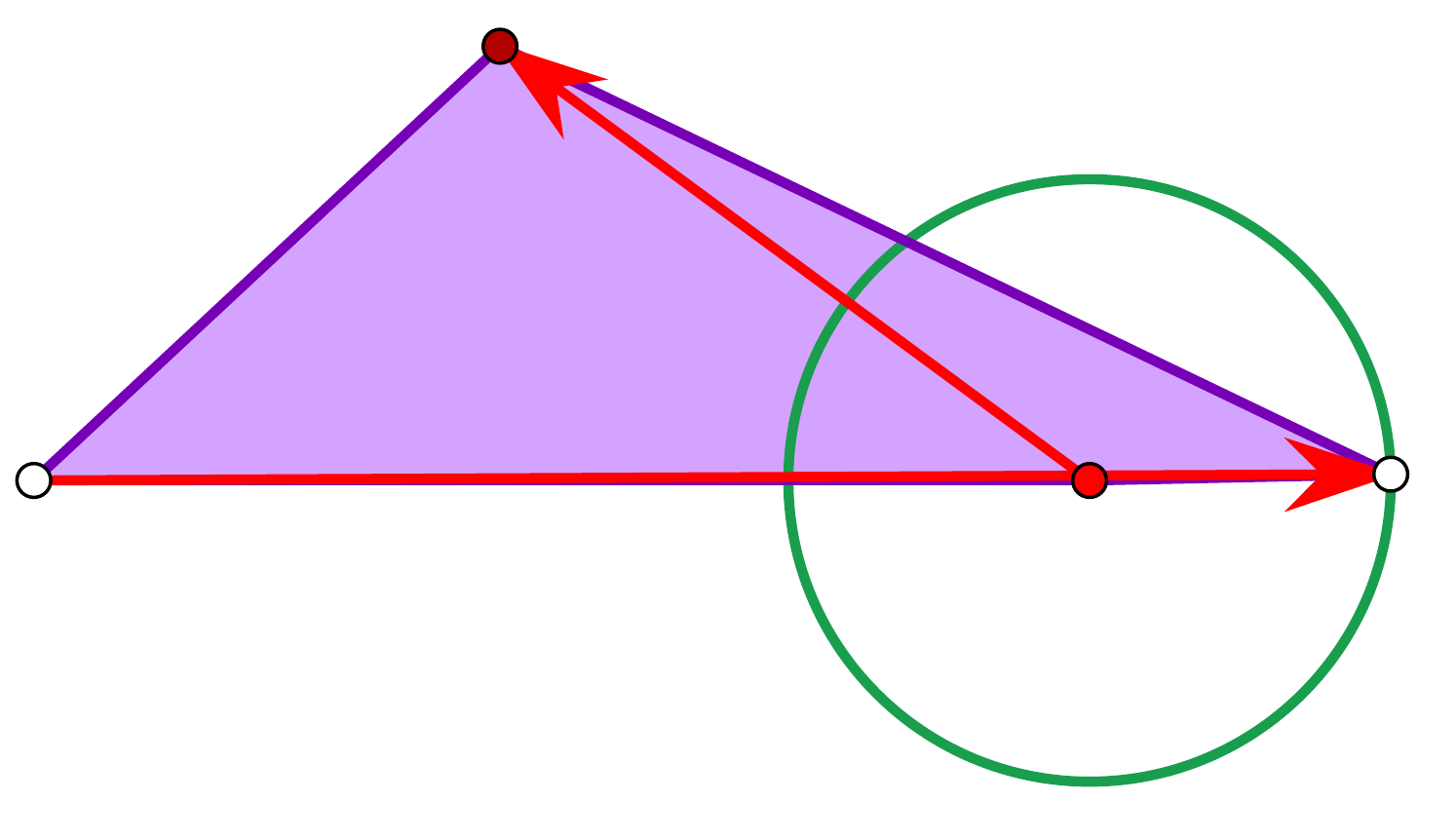}}
\hspace{3pt}
{\includegraphics[width=0.22\textwidth]{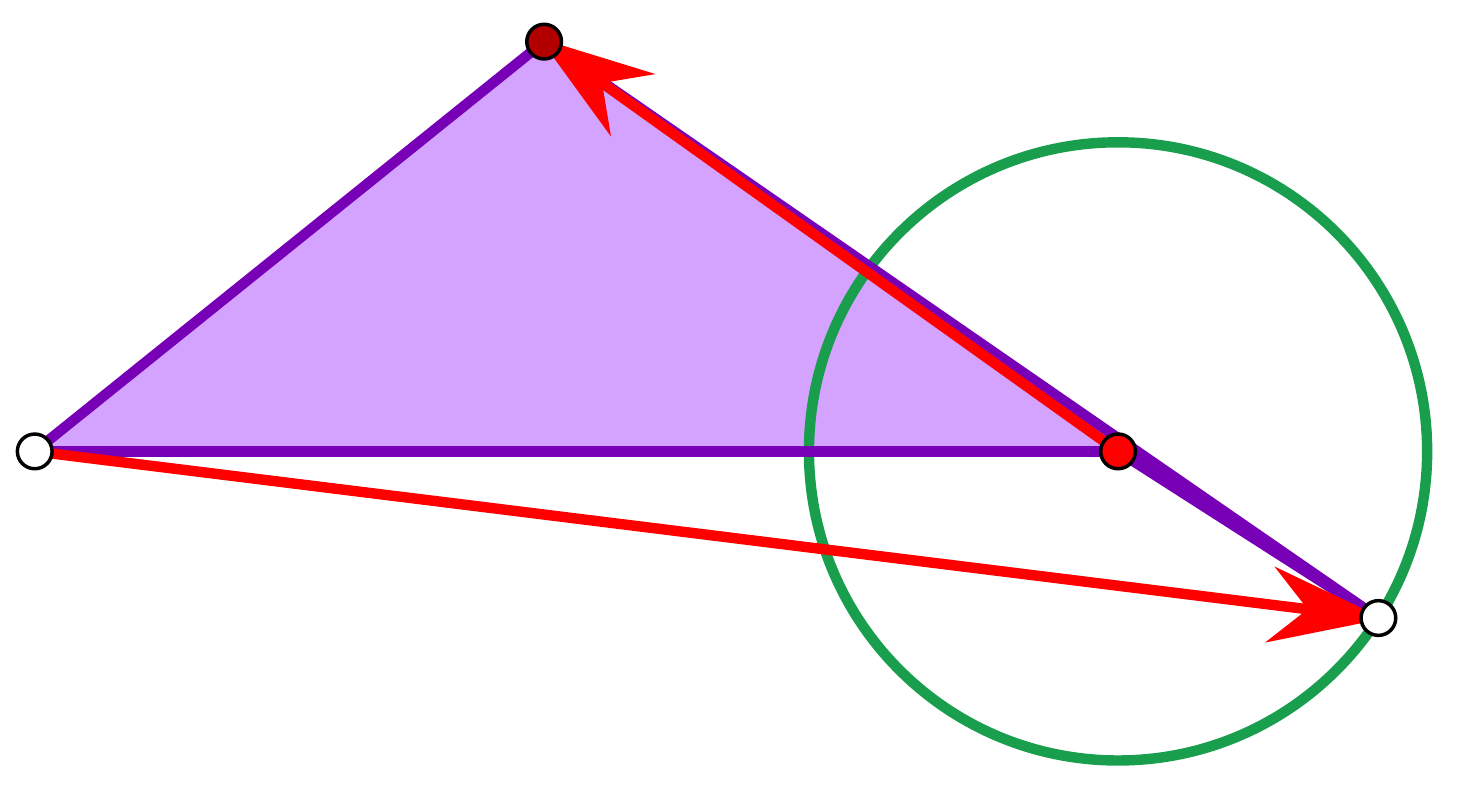}}
\caption{Auxetic-to-non-auxetic transition configurations: as the second edge rotates clockwise around the fixed first edge, there is an auxetic interval between the two configurations on the left and another one (in reverse) between the two configurations on the right.}
 \label{fig:Transitions}        
 \end{figure}

\noindent
The two theorems are illustrated in Fig.~\ref{fig:4barQ}.
The four transition configurations (from auxetic to non-auxetic behavior) of the quadrilateral are shown in Fig.~\ref{fig:Transitions}.

\section{The elliptope and other spectrahedra}
\label{sec:spectrahedra}

When aiming at more general results, an immediate problem is that of an adequate 
description of the configuration space of a given periodic framework. For a couple of instances where global configuration spaces have been determined, we refer to \cite{BS1,BScccg,BS8}. In the present context, as we focus on periodic frameworks with two vertex orbits, we rely on a procedure initiated in \cite{BS8}.

\medskip
In arbitrary dimension $d$, a connected periodic framework with two vertex orbits and all edges with endpoints in different orbits is completely determined by a vertex and the edge vectors emanating from it. The endings of these vectors are all in the opposite orbit and  the periodicity lattice is generated by vectors between pairs of such endings. We fix the first vertex at the origin and denote the edge vectors
emanating from it by $v_i$, $i=0,1,...,m-1$. In dimension $d$, connectedness requires
$m\geq d+1$. We may and shall assume that the labeling is such that the period vectors $\mu_j=v_j-v_0$ form a basis of ${R}^d$ for $j=1,...,d$. Then, period vectors $\mu_k=v_{k+1}-v_1$, $k > d$ have unique expressions with rational coefficients $a_{ki}$:

\begin{equation}\label{Q-dependence}
\mu_k=\sum_{i=1}^{d} a_{ki}\mu_i, \ \ \ k > d
\end{equation}

\noindent
These relations must be preserved in any periodic deformation of the framework and
together with the first $d+1$ vectors $v_i$ will determine the framework.

\medskip
\noindent
Up to isometry, the configuration of vectors $v_i\in {R}^d$, $i=0,1,...,d$ is
determined by their Gram matrix $G=V^tV$, where $V$ is the $d\times (d+1)$ matrix with column vectors $v_i$. This $(d+1)\times (d+1)$ symmetric matrix is positive semidefinite of rank $d$. If we take into account the prescribed squared lengths of the edges, say: 

\begin{equation}\label{sql}
\langle v_i, v_i \rangle = s_i > 0
\end{equation}

\noindent
we see that, as the framework deforms, the corresponding matrix $G=(g_{ij})$ 
belongs to the affine section of the positive semidefinite cone defined by
$g_{ii}=s_i$, $i=0,1,...,d$. This is a {\em spectrahedron} linearly equivalent
whith the spectrahedron of {\em correlation matrices}, also known as the {\em elliptope}, which corresponds to all diagonal entries equal to one \cite{LP,V}.
More precisely, $G$ belongs to the boundary of this spectrahedron, which is part
of the algebraic hypersurface of degree $d+1$ given by the vanishing of the determinant. Now, we use the fact that relations (\ref{Q-dependence}) and (\ref{sql}) for $k > d$ amount to affine sections of the spectrahedron and obtain:

\begin{theorem}\label{param}
The configuration space of a connected periodic framework in ${R}^d$ with
two vertex orbits has a natural compactification as the boundary of 
a spectrahedron, which is itself linearly equivalent with an affine section
of an elliptope. Topologically, this compactification is a sphere. 
\end{theorem}

\noindent
This generalizes the result in section~\ref{sec:4bar}, where the configuration space was a circle or a circle without a point.

\section{Equations in lattice coordinates}
\label{sec:2vertex}

Lattice coordinates are linear coordinates relative to a basis made of independent
generators for the periodicity lattice \cite{BS3,BS10}. In crystallography they are called {\em fractional coordinates}. Suppose, as in Definition~\ref{A-path}, that
the chosen periodicity generators form the $d\times d$ matrix $\Lambda$, with
associated Gram matrix $\omega=\Lambda^t \Lambda$. Metric relations in lattice coordinates are expressed via $\omega$. For periodic frameworks with two vertex orbits, one orbit will be represented in lattice coordinates by the integer lattice ${Z}^d$ and the other one by $q+{Z}^d$, for some $q\in {R}^d$. We have $(q,\omega)\in {R}^d\times {R}^{d(d+1)/2}$.

\medskip \noindent
Relating to notations used in the previous sections, the period vectors $\mu_i$ 
become vectors $n_i\in {Z}^d$ and the vertex used there as origin has now lattice
coordinates $q$. Thus, the edge vectors emanating from this vertex are $n_i-q$,
$i=0,1,...,m-1$. The system of equations expressing the constancy of the squared length
of the edges becomes:

\begin{equation}\label{syst}						
\langle \omega (n_i-q), n_i -q \rangle = s_i
\end{equation}

\noindent
This gives:

\begin{equation}\label{syst1}						
\langle \omega n_i, n_i \rangle -2 \langle \omega q, n_i \rangle + 
\langle \omega q, q \rangle= s_i
\end{equation}

\medskip \noindent
and leads to the equivalent system:

$$ \langle \omega q, q \rangle=s_0 $$
\begin{equation}\label{syst2}						
 2 \langle \omega q, n_i \rangle = 
  s_0 - s_i + \langle \omega n_i, n_i \rangle, \ \ i=1,...,d
\end{equation}
$$ 2 \langle \omega q, n_j \rangle =  s_0 - s_j + \langle \omega n_j, n_j \rangle, \ \ j > d.  $$

\medskip \noindent
Expressing $n_j, j> d$, as linear combinations (with rational coefficients) of
$n_i, i=1,...,d$ via (\ref{Q-dependence}), equations (\ref{syst2}) can be used to eliminate $\omega q$
from the subsequent equations which become linear constraints just for $\omega$.
 
\medskip \noindent
Moreover, the system (\ref{syst2}) can be used to write all coordinates of 
$\omega q$ as degree one expressions in the coordinates of $\omega$. Substituting in the first equation and rewriting (\ref{syst2}) as:

\begin{equation}\label{syst3}						
 2 \langle q,  \omega n_i \rangle = 
  s_0 - s_i + \langle \omega n_i, n_i \rangle, \ \ i=1,...,d
\end{equation}

\noindent
we obtain a system of
$d+1$ equations in $(q,\omega)$, where $\omega$ is already restricted to the
affine subspace described above. The equations have degree one in $q$ as well as in $\omega$.

\medskip \noindent
If we return to the case $d=2$, $m=4$ of section~\ref{sec:4bar}, we have a curve defined by three equations in $(q,\omega)\in {R}^2\times {R}^2$. Projection on
$q\in {R}^2$ gives a planar cubic curve and allows a graphical recognition of 
the auxetic intervals based on Theorem~\ref{ellipse}.

\section{Conclusion}
\label{sec:conclusion} 
 
In this paper we studied the problem of finding auxetic regions in the global deformation space of a periodic framework. Points in such regions allow local auxetic one-parameter deformations. We focused on periodic frameworks with two vertex orbits in arbitrary dimension $d$, for which the case of planar frameworks associated to four-bar mechanisms served as an intuitive example. We described algebraic techniques and geometric perspectives for addressing the problem.

\vspace{0.3in}
\noindent
{\bf Acknowledgement}
 \ This work was supported by the National Science Foundation (awards no. 1319389 and 1704285 to C.S.B.,\  awards no. 1319366 and 1703765 to I.S.)\ \ and the  National Institutes of Health (award 1R01GM109456 to C.S.B. and I.S.).

\vspace{0.4in} \noindent
\textsc{Ciprian S. Borcea}, Department of Mathematics, Rider University, Lawrenceville, NJ 08648, USA  

\medskip \noindent 
 \textsc{Ileana Streinu}, Computer Science Department, Smith College, Northampton, MA 01063, USA.

\end{document}